\documentclass[a4paper,12pt]{article}
\usepackage[latin1]{inputenc}
\usepackage[T1]{fontenc}
\usepackage[english]{babel}
\usepackage{amssymb}
\usepackage{latexsym}
\usepackage{xcolor}
\usepackage{amsmath}
\usepackage{amsthm}
\usepackage{mathrsfs}
\usepackage{amsxtra}
\usepackage{eucal}
\usepackage[left=1.5cm,right=1.5cm,top=1cm,bottom=1.5cm]{geometry}
\usepackage{tikz}
\usetikzlibrary{patterns}
\usepackage{amsfonts}
\usepackage{lmodern}
\usepackage{graphicx}
\usepackage{float}
\providecommand{\LyX}{L\kern-.1667em\lower.25em\hbox{Y}\kern-.125emX\@}

\def\R{\mathbb{R}}
\def\E{\mathbb{E}}
\def\dint{\displaystyle\int}

\def\P{\mathbb{P}}
\newtheorem{thm}{Theorem}[section]

\newtheorem{claim}[thm]{Claim}

\newtheorem{cor}[thm]{Corollary}

\newtheorem{lem}[thm]{Lemma}

\newtheorem{prop}[thm]{Proposition}
\newtheorem{rk}[thm]{Remark}

\newenvironment{class}[1][AMS Classification]{\textbf{#1.} }{}
\newenvironment{MC}[1][Key words]{\textbf{#1.} }{}

\def\disp{\displaystyle{}}

\def\al{\alpha}
\def\si{\sigma}
\def\be{\beta}

\title{\textbf{  On Malliavin differentiability and absolute continuity of one-dimensional doubly perturbed diffusion processes}}

\author{ R. BELFADLI  $^{1}$ $\quad$ L. BOULANBA   $^{2}\quad$ Y. OUKNINE $^{3}$ $^{4}\quad$
\vspace*{0.1in}\\
$^{1}$  Department of Mathematics, Faculty
of Sciences and Technologies,\\ Cadi Ayyad University, Marrakesh, Morocco.\\
rachid.belfadli@uca.ac.ma 
\vspace*{0.1in}\\
$^{2}$ Department of Mathematics, CRMEF SM.\\ Laboratory  LISTI, ENSA-Agadir,  Ibn Zohr University, Morocco.\\
\quad l.boulanba@crmefsm.ac.ma 
\vspace*{0.1in}\\
$^{3}$  Department of Mathematics, Faculty of Sciences Semlalia.\\
Cadi Ayyad University, Marrakesh, Morocco. \\
ouknine@uca.ac.ma
\vspace*{0.1in}\\
$^{4}$ Mohammed VI Polytechnic University, Africa Business School,\\
Madinat  Al Irfane, RABAT, Morocco. \\
 youssef.ouknine@um6p.ma
}

\begin{document}
\maketitle \textbf{Abstract.} In this paper, we establish Malliavin differentiability and absolute continuity for  $\alpha, \beta$-doubly perturbed diffusion  process with parameters  $\alpha <1$ and $\beta <1$   such that $|\rho|  < 1$, where $ \rho : = \frac{\alpha\beta}{(1-\alpha)(1-\beta)}$. Furthermore, under some regularity conditions on the coefficients, we prove that the solution $X_t$ has a smooth density for all $t\in(0, t_0)$ for some finite number $t_0>0$.  Our results recover earlier works by Yue and Zhang \cite{Yue2015} and Xue, Yue and Zhang \cite{Xue2016},  and  the proofs are based on the techniques of the Malliavin calculus.
\\

\begin{MC}
Doubly perturbed stochastic differential equations; Malliavin calculus; Absolute continuity; Smooth density.
\end{MC}

\begin{class} 60H07;  60H10.
\end{class}
\section{Introduction}
 Stochastic differential equations (SDEs) involving the past maximum and/or the past minimum processes has attracted  much interest and has been  treated by several authors during the past few decades, (see for instance \cite{LGY1990, CPY95, Davis1996, PW1997, CD1999, CDH2000, Doney-Zhang-2005, Belfadli2009, MHM2018, GM2023, Belfadli2024}).
 
The present work is a continuation of studies concerned with perturbed stochastic diffusion processes. Namely,  let $\Omega= C(\mathbb{R}^{+}; \mathbb{R})$ be the set of continuous functions from $\mathbb{R}^{+}$ into $\mathbb{R}$,
$\P$ the Wiener measure on $\Omega$, $(W_t)_{t\geq 0}$ the process of
coordinate maps from $\Omega$ into $\mathbb{R}$, ${\cal
F}=\sigma\{W_t,\,t\geq 0\}$, $({\cal F}_t)_{t\geq 0}$ the completion
of the natural filtration of $W$ with the $\P$-null sets of $\cal F$.
Therefore $(W_t)_{\geq 0}$ is a standard Brownian motion on the
filtered probability space $(\Omega, {\cal F}, ({\cal F}_t)_{t\geq
0},\P)$. We consider the following one-dimensional doubly perturbed diffusion process
 \begin{eqnarray}\label{e2}
X_t=x+\int_0^t \si(X_s)dW_s +\int_0^t b(X_s)ds +\al\max_{0\leq
s\leq t}X_s +\be \min_{0\leq s \leq t}X_s, \, \, \, t\geq 0,
\end{eqnarray} 
where $b$ and $\sigma$ are real valued Lipschitz functions. 

Let us first briefly comment on Eq. (\ref{e2}) and the  permissible values of  parameters $\alpha$ and $\beta$. Recall that, under the conditions $\alpha <1$, $\beta <1$ and $|\rho|< 1$, existence and uniqueness of a strong solution for the above equation was established by Carmona, Petit and Yor \cite{CPY95} in the case $x=0, \sigma \equiv 1, b\equiv 0$, and by Belfadli, Hamad\`{e}ne and Ouknine \cite{Belfadli2009} with Lipschitz coefficients $b$ and $\sigma$. Studying the same functional equation as in \cite{CPY95}, Davis \cite{Davis1996} proved that the result still holds in the remaining case when $|\rho| = 1$. Moreover, it is proved in \cite{CPY95} that no solution exists if $\alpha \geq 1$ or $\beta \geq 1$. In the same direction, and with non Lipschitz coefficients, a version of Eq.  (\ref{e2}) involving a jump term was studied by Luo \cite{Luo2009}. The author obtained existence and uniquenss of the solution under a more restrictive condition on  $\alpha$ and $\beta$, namely   $|\alpha| + |\beta| < 1$ ( see condition $(A5)$ in \cite{Luo2009}).
 A multidimensional extension  of (\ref{e2}) has been recently considered by Gassiat and M\k{a}dry \cite{GM2023}.  For fixed integer $d\geq 1$ and assuming  $ \alpha_i<1, \  \beta_i<1, \ |\rho_i|<1$, for all $i=1, \cdots, d$, the authors established existence, pathwise uniqueness and stability estimates for a $d$-dimensional version of (\ref{e2}).  

\begin{figure}[H]
\centering

\begin{tikzpicture}[scale=1.5]

 \draw[thick,->,line width=0.5mm] (0,0) -- (0,2) node[left] {$\beta$};
 \draw[thick,->,line width=0.5mm] (0,0) -- (2,0) node[below ] {$\alpha$}; 
 \draw[line width=0.5mm] (0,0) -- (-2,0);  
 \draw[line width=0.5mm] (0,0) -- (0,-2); 
  \draw[line width=0.2mm,dashed] (0,0.5) -- (-2,0.5);  
 \draw[line width=0.2mm,dashed] (0.5,0) -- (0.5,-2);
  \draw[line width=0.2mm,dashed] (1,0) -- (0,-1); 
   \draw[line width=0.2mm,dashed] (0,1) -- (-1,0);
   \draw[line width=0.2mm,dashed] (-1,0) -- (0,-1);  
 \foreach \x in{1,-1 }
 \draw[line width=0.3mm] (\x , 1pt) -- (\x , -2pt) node[anchor=north] {$\x$};
 
 \foreach \y in{1, -1}
 \draw[line width=0.4mm] (2pt,\y) -- (-2pt,\y) node[anchor=east] {$\y$};
\node[left] at(0,0.45) {$ \frac{1}{2} $}; 
\node[below] at (0.6,0) {$ \frac{1}{2} $}; 
%\node at (0.61,0.61){\rotatebox{-45}{$\alpha+\beta -1=0$}};
\node at (-1.5,0.83){$\beta=\frac{\alpha-1}{2\alpha -1}$};
\node at (1.3,-1.5){$\beta=\frac{\alpha-1}{2\alpha -1}$};
 \draw [line width = 1.2pt, domain=-2:0, samples=100,>=latex]%
plot (\x, {(1-\x)/(1-2*\x)}) ;
\draw [line width = 1.2pt, domain=0.6:1, samples=100,>=latex]%
plot (\x, {(1-\x)/(1-2*\x)}) ; 
 
 \draw [line width = 1.2pt, domain=0:1, samples=100,>=latex]%
plot (\x, {1-\x}) ; 
\draw[pattern = north east lines,opacity=.5, line width = 1.2pt,draw=none] (-2,0.533)%
plot[domain=-2:0] (\x, {(1-\x)/(1-2*\x)})--(0,-2)--(-2,-2)--cycle;
\draw[pattern = north east lines,opacity=.5, line width = 1.2pt,draw=none] (0.5,-2)%
plot[domain=0.6:1] (\x, {(1-\x)/(1-2*\x)})--(1,0)--(0.5,0.5)--cycle;
\draw[pattern = north east lines,opacity=.5, line width = 1.2pt,draw=none] (0,1)%
plot[domain=0:0.6] (\x, {(1-\x)})--(0.5,0.5)--(0.6,-2)--(0,-2)--cycle;

\end{tikzpicture}

\caption{Domain of valid shape parameters $\alpha$ and $\beta$} \label{figure1}
\end{figure}
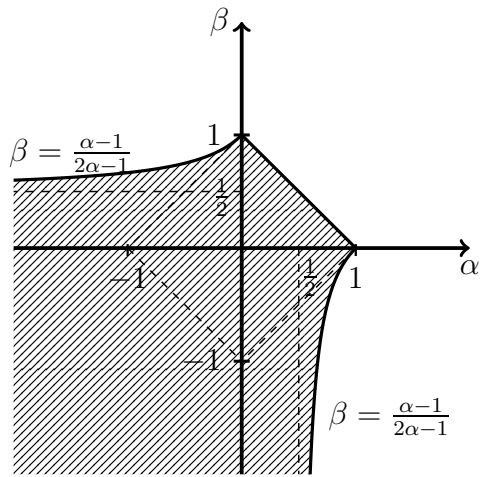

Henceforth, with these  comments in mind, it seems that the natural conditions to consider in every aspect of study of Eq. (\ref{e2})  are $\alpha<1$, $\beta <1$ and $|\rho|<1$. For convenience of the reader we illustrate in Figure \ref{figure1} below the domain of valid parameters $\alpha$ and $\beta$.  It should be noted that the results of the current paper are obtained in that region of parameters, and this constitute one of the more challenging contribution of this paper.
 
Next,  it is worth nothing that when $\beta=0$; that is in the  singly perturbed case and when the SDE involve only the maximum process, absolute continuity of the laws of solutions was established by Yue and Zhang in  \cite{Yue2015}, and  smoothness of the density  was addressed by Xue, Yue and Zhang in  \cite{Xue2016}.  In the same spirit of results,  the purposes of this paper are threefold. Firstly, to prove Malliavin differentiability of the solution; namely, $X_t \in \mathbb{D}^{1,2}$, for each $t > 0$. Secondly, to establish the absolute continuity of law of the solutions with respect to the Lebesgue measure, and lastly to proving the smoothness of the density of the solution. This is achieved mainly by using techniques of Malliavin calculus.

Our results cover and improve those established in \cite{Yue2015} and \cite{Xue2016} corresponding to the singly perturbed case.  More importantly, our proofs are more challenging. In fact,  the core  of proofs requires some sophisticated estimates which rely heavily on applying Skorokhod's lemma and Picard approximation. This is in fact the main technical part of the paper. 

The content of this paper is organised as follows. Section 2 contains preliminary notions and results underpinning our proofs, while Section 3 is  mainly devoted to prove the Malliavin differentiability of $X_t$ for any $t\geq 0$. Regarding absolute continuity of the laws of the solutions, it constitutes the subject of Section 4. Finally, in Section 5, we study the smoothness of the density for the solution of the considered equation.

Throughout this paper, $c$, $c^{\prime}$, $c^{\prime \prime}, \ldots$ will denote constant  whose value may vary from one expression to the other and $|| f||_{\infty}= \sup_{x}|f(x)|$ the norm of a bounded function $f$.

\section{Preliminaries} 
In this section we present some preliminaries on Malliavin calculus and recall some results which  will play important role in the analysis of regularity and absolute continuity of the laws for the solutions to (1).  For a more complete exposition, we refer to \cite{Nualart2006}. 

Let $H$ be the Hilbert space $\mathrm{L}^{2}(\mathbb{R}^{+};\mathbb{R})$, and for each $h\in H$ and we set 
$$
W(h) : = \int_{0}^{+\infty}h(s)dW(s).
$$

\noindent The isonormal Gaussian process $\{W(h), \ h\in H\}$  plays a crucial role in the definition of the derivative operator.  Let $\mathscr {S}$ be the set of all smooth cylindrical random variables of the form
\[F=f(W(h_{1}), W(h_{2}),\ldots, W(h_{n}))\]
where $n\geq 0$, $f\in C^{\infty}$ has a compact support and $h_{i}\in H$. The Malliavin derivative of $F$ with respect to $W$ is the element belonging to $L^{2}(\Omega; H)$ defined by
\[D_sF=\sum_{i=1}^{n}\frac{\partial f}{\partial x_{i}}(W(h_{1}), W(h_{2}),\ldots, W(h_{n}))h_{i}(s), \quad s\geq 0.\]

This operator can be extended to the closure $\mathbb{D}^{1,2}$ of $\mathscr S$ with respect to the norm
\[\Vert F\Vert^{2}_{1,2}:= \E[F^{2}] + \E[||D_{\cdot}F||_{H}^{2}].\]
As usual, for  $k \geq 2$, one can define by iteration the $k$th derivative $D^kF$, which is an element of $L^{2}(\Omega; H^{\odot k})$. The domain $\mathbb{D}^{k, 2}$ denotes the closure of $\mathscr{S}$ with respect to the norm  $||\cdot ||_{k, 2}$ defined by
\[\Vert F\Vert^{2}_{k,2}:= \E[F^{2}] + \sum_{i=1}^{k}\E[||D^{i}_{\cdot}F||_{H^{\otimes i}}^{2}].\]
See \cite{Nualart2006} for a more general definition of higher order Malliavin derivatives and spaces $\mathbb{D}^{k,p}$. 

We will make use of the following crucial result in the sequel, which is a slight modification of Lemma $1.2.3$ in \cite{Nualart2006} for higher order $k\geq 1$.
\begin{prop}\label{appro:D12} Let $k \geq 1$ be an integer.
Let $\{F_n, n \geq 1\}$ be a sequence of random variables in $\mathbb{D}^{k,2}$ that converges to $F$ in $\mathrm{L}^{2}(\Omega)$ and such that
$$
\sup_{n}\E[\|D^{k}F_n\|^2_{H^{\otimes k}}]< \infty.
$$
Then, $F$ belongs to $\mathbb{D}^{k,2}$, and the sequence of derivatives $\{DF_n, n \geq 1\}$ converges to DF in the weak topology of $\mathrm{L}^{2}(\Omega; {H^{\otimes k}}).$

\end{prop}
 The Malliavin derivative satisfies the following chain-rule. If $F=(F_{1},\ldots,F_{n})$ is a random vector with components in $\mathbb{D}^{1,2}$ and if $\varphi :\R^{n}\rightarrow \R$ is a continuously differentiable function with bounded partial derivatives, then
$\varphi(F_{1},...,F_{n})\in \mathbb{D}^{1,2}$ and we have, for any $s\geq 0$:
\begin{equation} \label{chaine rule:e1}
D_s\varphi(F_{1},...,F_{n})=\sum_{i=1}^{n}\frac{\partial \varphi}{\partial x_{i}}(F_{1},\ldots,F_{n})D_sF_{i}.
\end{equation} 
The chain-rule (\ref{chaine rule:e1}) still holds when $\varphi$ is only Lipschitz continuous with Lipschitz constant $K$. Indeed, there exists a random vector $Y=(Y_1, \ldots, Y_n)$ with $|Y|\leq K$ almost surely, such that
\begin{equation} \label{chaine rule:e2}
D\varphi(F_{1},...,F_{n})=\sum_{i=1}^{n}Y_iDF_{i}.
\end{equation} 
%We denote by $\delta$ the divergence operator; that is the adjoint of the operator $D$. A random element 

\noindent Sufficient conditions for the differentiability of the supremum of a continuous process are given by 
\begin{prop}\label{prop1:Nualart} (Proposition 2.1.10 in \cite{Nualart2006})
\noindent Fix  $T>0$ and let  $X = \{X_t, \, 0\leq  t \leq T\}$ be a continuous process. Suppose that 
\begin{itemize}
\item[(i)] $\E[\sup_{t\in [0, T]} X_t^{2}] < \infty;$
\item[(ii)] for any $t \in [0, T]$, $X_t \in \mathbb{D}^{1,2}$, the $H$-valued process $\{DX_t, 0\leq t \leq T\}$ 
possesses a continuous version, and $\E[\sup_{t\in [0, T]} \|DX_t\|_{H}^{2}] < \infty.$
\end{itemize}
Then, the random variable $M := \sup_{t\in [0, T]} X_t$ belongs to $\mathbb{D}^{1,2}$, and we have a.e.
\begin{equation}\label{nom_max=e1}
||  D( \sup_{0 \leq t \leq  T} X_t) ||_H^2 \leq  \max_{0\leq t \leq T}||DX_t ||_H^2.
\end{equation}
\end{prop}

The following result on the absolute continuity of the law of a one-dimensional random variable is a well-known fact of the Malliavin calculus
\begin{prop} \label{absl}
Let $F\in\mathbb{D}^{1,2}$. If $\|DF\|_{H} > 0$ a.s., then the law of $F$ is absolutely continuous with
respect to the Lebesgue measure on $\mathbb{R}$.
\end{prop}

  A general criterion for smoothness of densities for random vectors belonging to $\mathbb{D}^{\infty} : = \bigcap_{p, k\geq 1}\mathbb{D}^{k,p}$ is given by

\begin{prop}(Proposition 2.1.5 in \cite{Nualart2006})
Let $F = (F_1, \cdots, F_m)$ be a random vector satisfying the following conditions : 
\begin{itemize}
\item[$i)$] $F^{i}$ belongs to $\mathbb{D}^{\infty} $ for all $i=1, \cdots, m$;
\item[$ii)$] The Malliavin matrix $\gamma_{F} = (\langle F_{i} , F_{j}\rangle_{H})_{1\leq i, j \leq m}$ satisfies 
$$
(\text{det}\  \gamma_{F})^{-1} \in \bigcap_{p\geq 1} {L}^{p}(\Omega).  
$$
\end{itemize}
Then  $F$  possesses an infinitely differentiable density on $\mathbb{R}^m$.
\end{prop}  
The following result ensues for one-dimensional random variables
\begin{cor}\label{cor:1} If $F$ is a random variable belonging to $\mathbb{D}^{\infty}$ such that $||DF||_{H}\geq c$ a.s., for some positive constant $c$, then $F$  possesses an infinitely differentiable density on $\mathbb{R}$.
\end{cor}  
%{{\color{blue}{
\section{Malliavin differentiability of $\alpha, \beta$-doubly perturbed diffusion processes }
Let us consider, for $\alpha <1$ and $\beta<1$ such that $|\rho|<1$, the solution of the following $\alpha, \beta$-doubly perturbed SDEs :
\begin{eqnarray}\label{pe2}
X_t=x+\int_0^t \si(X_s)dW_s +\int_0^t b(X_s)ds +\al\max_{0\leq
s\leq t}X_s +\be \min_{0\leq s \leq t}X_s,
\end{eqnarray}
where from this point forward, we will assume $b$ and $\sigma$ are two real valued Lipschitz functions.  The  main result of this section can be stated as follows
\begin{thm} \label{thm1} Let $\alpha <1$ and $\beta<1$ such that $|\rho|<1$. 
Let $(X_t)_{t\geq 0}$ be the unique solution of Eq. (\ref{e2}). Then,  $X_t $ belongs to $\mathbb{D}^{1, 2}$ for any $t\geq 0$. 
\end{thm}

\begin{rk}This result cover and  improve the one studied by Yue and Zhang \cite{Yue2015} in which $\beta$ is just $0$.
\end{rk}

The main idea of the proof  is essentially based on Proposition \ref{appro:D12} and use Picard iterations scheme.  More precisely, we consider the sequence of processes implicitly defined by : 
\begin{equation}\label{eq2}
\left\{
\begin{array}{l}
  X_t^0 = \frac{x}{1-\alpha},\\
X_t^{n+1} = x +  \int_0^t \sigma(X_s^n)dW_s  + \int_0^tb(X_s^n)ds + \alpha.M_t^{n+1} + \beta.I_t^{n+1},  
\end{array}
\right.
\end{equation}
where $M_t^{n+1} := \max_{0\leq s \leq t} X_{s}^{n+1}$  and   $I_t^{n+1} := \min_{0\leq s \leq t} X_{s}^{n+1}$,  for every $t\geq 0$ and $n \geq 0$.

\vspace{0.2cm}
\noindent As it set forth in \cite{Belfadli2009} by relying on Skorokhod's lemma the sequence $(X^n)_{n\geq 0}$ is well defined and converges uniformly on compact intervals almost surely to $X$. Furthermore, 
\begin{equation}\label{convergence1}
\lim_{n\longrightarrow \infty}\mathbb{E}[|X_t^n - X_t|^2] = 0, \ \text{for all} \  t\geq 0. 
\end{equation}
\subsection{Auxiliary results}
As a preliminary of the proof of Theorem \ref{thm1}, we need the following 

 \begin{prop}\label{proposition1} Let $(X^n)$ be given by (\ref{eq2}). Then, for all integer $n \geq 0$,
 
 \begin{equation}\label{eq3}
X_t^{n}\in \mathbb{D}^{1, 2}, \ \text{and} \  \mathbb{E}\left[\int_{0}^{t}\|D X_s^n\|_H^{2}ds \right]< \infty,\,\, \mbox{for all}\, \, t\geq 0. 
\end{equation} 
 \end{prop}
 \noindent {\bf Proof of Proposition \ref{proposition1}}: We employ  an induction argument on $n$. Obviously, (\ref{eq3}) holds for $n=0$. Let $n\in \mathbb{N}$. Assuming  $(\ref{eq3})$ to hold for $n$, we will prove it for $n+1$. Applying the chain-rule (\ref{chaine rule:e2}) to $X_t^n$ and to the functions $\sigma$ and $b$, we deduce that $\sigma(X_t^n)$  and $b(X_t^n)$ belong to $\mathbb{D}^{1, 2}$, and there exist adapted processes $\bar{\sigma}^n(\cdot)$ and  $\bar{b}^n(\cdot)$ uniformly bounded by $K$, such that   
\begin{equation}\label{eq44}
D_{r}\sigma(X_s^n) = \bar{\sigma}^n(s)D_{r}(X_s^n)1_{r\leq s},
\end{equation} 
\begin{equation}\label{eq5}
D_{r} b(X_s^n) = \bar{b}^n(s)D_{r}(X_s^n)1_{r\leq s}.
\end{equation}  

\noindent Thus, the processes $\{D_r\sigma(X_s^n),  s\geq r \}$ and $\{D_r b(X_s^n),  s\geq r \}$ are square integrables and adapted, and from (\ref{eq44}) and (\ref{eq5}) we get   
\begin{equation}
|D_{r}\sigma(X_s^n)| \leq K |D_{r} X_s^n|,
\end{equation} 

\begin{equation}\label{derivation:b}
|D_{r}b(X_s^n)| \leq K |D_{r} X_s^n|.
\end{equation}

Consequently, by Lemma 1.3.4 and Proposition 1.3.8 in \cite{Nualart2006}, we obtain that 
\begin{equation}
\int_{0}^{t}\sigma(X_s^n)dW_s \in \mathbb{D}^{1,2}, \ \text{and} \ D_r\left(\int_{0}^{t}\sigma(X_s^n)dW_s\right) =  \sigma(X_r^n)1_{r\leq t} + \int_{r}^{t}D_{r}(\sigma(X_s^n))dW_s.
\end{equation}
Similarly, we have
\begin{equation}
\int_{0}^{t}b(X_s^n)ds \in \mathbb{D}^{1,2}, \ \text{and} \ D_r\left(\int_{0}^{t}b(X_s^n)ds\right) =  \int_{r}^{t}D_{r}(b(X_s^n))ds.
\end{equation}
Therefore 
\begin{equation}\label{eq:a reg}
a^n_t:=x+ \int_{0}^t \sigma(X^n_s) dW_s +  \int_{0}^t b(X^n_s) ds \in  \mathbb{D}^{1,2},
\end{equation}
and 
\begin{equation}\label{eq:a}
D_r(a^n_t)= \sigma(X_r^n)1_{r\leq t} + \int_{r}^{t}D_{r}(\sigma(X_s^n))dW_s + \int_{r}^{t}D_{r}(b(X_s^n))ds.
\end{equation}
In view of Eq. (\ref{eq2}), and in order to prove that $X^{n+1}_t \in \mathbb{D}^{1,2}$, it remains to establish that $M^{n+1}_t$ and $I^{n+1}_t$ belong to $\mathbb{D}^{1,2}$. 
For sake of comparison, let us first  make the following Remark.
\begin{rk}  In \cite{Yue2015}; that is when $\beta =0$, the employed  proof  consist in expressing $M^{n+1}_t$ as the maximum of the process $\{\frac{a^{n}_s}{1- \alpha}, s \geq 0\}$ over the time interval $[0, t]$ and then  applying Proposition \ref{prop1:Nualart}. However, in our setting one cannot follows  the same pattern as in \cite{Yue2015}. The snag is that neither $M^{n+1}_t$ nor $I^{n+1}_t$ can be expressed explicitly as  a functional of $\{X^{n}_s, s \leq t\}$. Indeed, $M_t^{n+1}$ and $I_t^{n+1}$ appear as  solutions of functional equations involving the process $X^{n+1}$ itself (see for instance, formulas (\ref{eq4}) and (\ref{e6}) below). To handle the aforementioned difficulty  and based on the contraction property of Skorokhod map  we use again a Picard approximation and introduce a sequence of continuous  adapted processes $(M^{n+1,m})_m$  converging uniformly on compact intervals almost surely to $M^{n+1}$ and then use Proposition  \ref{appro:D12} to conclude that $M_t^{n+1} \in \mathbb{D}^{1, 2}$.
\end{rk}
We state the following Lemma whose  proof is postponed until the proof of Proposition \ref{proposition1} is finished.
 \begin{lem} \label{lem1} Let $t\geq 0$  and $a^n_t$ be given by  (\ref{eq:a reg}). Then, we have  
 \begin{enumerate}
 \item[(i)]  $\max_{0\leq s \leq t}(a^n_s) \in \mathbb{D}^{1, 2}$ and \[||D (\max_{0\leq s \leq t}(a^n_s))||^2_H \leq \max_{0\leq s \leq t}||Da^n_s||^{2}_H\, \, \mbox{ a.s.} \]
Moreover, 
\begin{align}\label{eq:max an}
 \E\left[\max_{0 \leq s \leq t}|| Da^{n}_s||^2_{H}\right]\leq  c\left(1 +\int_{0}^t  \E[||DX^{n}_{u}||_H^2]du\right).
\end{align}
 \item[(ii)] $M^{n+1}_t \in \mathbb{D}^{1, 2}$  \,  and\,   $\E(\max_{0\leq s \leq t} || D M^{n+1}_s ||^2_H)\leq c \E(\max_{0 \leq s \leq t} || D a^n_s ||^2_H)$.
 \item[(iii)] Similarly, $I^{n+1}_t \in \mathbb{D}^{1, 2}$ \,  and\,   $ \E[\max_{0 \leq s \leq t} || D I^{n+1}_s ||^2_H]\leq c \E[\max_{0 \leq s \leq t} || D a^n_s ||^2_H]$.
 \end{enumerate}
 \end{lem}

 \noindent As an immediate consequence of the points $(ii)$ and $(iii)$ of  Lemma \ref{lem1} and (\ref{eq:a reg}), we  obtain 
 \[X^{n+1}_t= a^n_t  +\alpha M^n_t  + \beta I^n_t \in \mathbb{D}^{1, 2}. \]
 For the second claim in (\ref{eq3}), we have
 
\begin{align*}
& \mathbb{E}\left[\int_{0}^{t}\|D X_s^{n+1}\|_H^{2}ds\right] \leq  3 \left(  \int_{0}^{t} \E[||Da^{n}_s||_H^2]ds + \alpha^2 \int_{0}^{t} \E[||DM^{n}_s||_H^2]ds + \beta^2 \int_{0}^{t} \E[||DI^{n}_s||_H^2]ds\right)\\
& \leq 3 t (1+ \alpha^2 + \beta^2) \E\left[\max_{0\leq s \leq t} || D a^n_s ||^2_H\right] \\
 & \leq c\left(1 +\int_{0}^t  \E[|DX^{n}_{u}||_H^2]du\right)< \infty,
\end{align*}
 where we have used Lemma \ref{lem1} in the last two inequalities. This finishes the proof of Proposition \ref{proposition1}. Let us now return to the proof of Lemme \ref{lem1}.
\vspace{.2cm}

\noindent {\bf Proof of Lemma  \ref{lem1}}: To prove the first claim $(i)$, we will use Proposition \ref{prop1:Nualart}. \\ On one hand, notice first that 
 \begin{align}\label{max:e1}
& \mathbb{E}\left[\max_{0\leq s \leq t}(a^n_s)^2\right] \leq  3\left( x^2+ \mathbb{E}\left[\max_{0\leq s \leq t}\left( \int_{0}^s \sigma(X^n_s) dW_s\right)^2\right] + \mathbb{E}\left[\max_{0\leq s \leq t} \left(\int_{0}^t b(X^n_s) ds\right)^2\right]\right)\nonumber\\
& \leq c \int_{0}^t \E\left(1+ |X^{n}_r|^2\right) dr \leq c  \sup_{n} \int_{0}^t \E\left(1+ |X^{n}_r|^2\right) dr=: A<\infty,
\end{align}
where we have used Burkholder-Davis-Gundy's inequality, Lipschitz's assumption on $b$
and $\sigma$ and the fact that $X^n$ converges to $X$ uniformly with respect to time parameter (see for instance \cite{Belfadli2009, CPY95}).\\ On the other hand, thanks to (\ref{eq:a}) we have 
\begin{align*}
& ||Da^{n}_s||^2_{H}= \int_{0}^{s} |D_ra^{n}_s|^2dr\\
& \leq 3 \left(  \int_{0}^s |\sigma(X^n_r)|^2 dr  + \int_{0}^s \left| \int_{r}^t D_r(\sigma(X^n_u))dW_u\right|^2 dr + \int_{0}^s \left| \int_{r}^s D_r(b(X^n_u))du \right|^2 dr\right).
\end{align*}
Then
\begin{align*}
&\E\left[\max_{0 \leq s \leq t}|| Da^{n}_s||^2_{H}\right] \leq   \E\left[\max_{0 \leq s \leq t}\int_{0}^{s} |D_ra^{n}_s|^2dr\right]\\
& \leq 3\left( \int_{0}^t\E[|\sigma(X^n_r)|^2] dr  +  \E\left[ \max_{0\leq s \leq t}\left(\int_{0}^s \left| \int_{r}^s D_r(\sigma(X^n_u))dW_u\right|^2 dr\right)\right]\right.\\
 & \left. +\E\left[\max_{0\leq s \leq t}\left( \int_{0}^s \left| \int_{r}^s D_r(b(X^n_u))du \right|^2 dr\right)\right]\right).
\end{align*}
Lipshitz condition on $\sigma$, Burkholder-Davis-Gundy inequality, (3)-(4) and arguing as in the proof of Proposition 3.3,  we get
\begin{align}\label{max:e2}
&\E[\max_{0 \leq s \leq t}  || Da^{n}_s||^2_{H}] \leq  c\left(A +\int_{0}^t \int_{r}^t \E[|D_rX^{n}_{u}|^2]dudr\right)\nonumber\\
& = c\left(A + \int_{0}^t  \E[||DX^{n}_s||^2_{H}]ds\right)<\infty,
\end{align}
which is exactly (\ref{eq:max an}). In addition, taking into account of (\ref{max:e1}) and using (\ref{nom_max=e1}) of Proposition \ref{prop1:Nualart}  concludes the  proof of the first claim $(i)$.\\
Let us now prove the second assertion $(ii)$.   As it has been mentioned above and  in light of (\ref{eq4}), we introduce the following Picard iteration scheme for the process $M^{n+1}$ by setting 
\begin{equation}\label{eq3:appro}
\left\{
\begin{array}{l}
  M_s^{n+1, 0} := x, \\
M_s^{n+1, m+1} :=\frac{1}{1-\alpha}\max_{0\leq u \leq s}\left\{a^n_u + \frac{\beta}{\beta -1} \max_{0\leq v\leq u }(-a^n_v-\alpha M^{n+1, m}_v)\right\}, 
\end{array}
\right.
\end{equation}
for all integer $m$ and $0\leq s\leq t$.
\begin{rk}
Observe that when $\beta=0$, $(M^{n+1, m})_{m\geq 1}$ is a stationary sequence and is in fact reduced to the  process $\{\frac{1}{1-\alpha}\max_{0\leq u \leq s}a^n_u\}$. Hence, in the case $\beta=0$ our proof is reduced to the one given in \cite{Yue2015}.
\end{rk}
Arguing as in Carmona, Petit and Yor \cite{CPY95}, it is easily seen that $(M^{n+1, m})_m$ converges uniformly on compact intervals almost surely to $M^{n+1}$ and that
\begin{equation}\label{conv:L2}
\lim_{m\longrightarrow \infty}\mathbb{E}(|M^{n+1, m}_t- M^{n+1}_t|^2 )= 0, \ \text{for all} \  t\geq 0.
\end{equation}
As an intermediate result, we start with the following Lemma.
\begin{lem}\label{lem2} For all integer $m\geq 0$, we have
\begin{itemize}
\item[(i)]  $\E[\max_{0 \leq s \leq t}  (M^{n+1, m}_s)^2]<\infty$;
\item[(ii)]  For all $t\geq 0$, $ M_t^{n+1, m} \in \mathbb{D}^{1, 2}, \E[\max_{0 \leq s \leq t}  ||DM^{n+1, m}_s||_H^2]<\infty$. 
\end{itemize}
Furthermore,\,\,  $\sup_m\E[ ||DM^{n+1, m}_t||_H^2]<\infty$, for all $t\geq 0$.
\end{lem}
\noindent {\bf Proof of Lemma \ref{lem2}}: By induction on $m$. The first assertion $(i)$ is obvious. In fact, using (\ref{eq3}) we have
\[\sup_{0\leq s \leq t}|M_s^{n+1, m+1}| \leq \left(\frac{1}{1-\alpha} + \frac{|\beta|}{(1-\alpha)(1-\beta)}\right)\sup_{0\leq s \leq t}|a_s^{n}| + \frac{|\alpha\beta|}{(1-\alpha)(1-\beta)}\sup_{0\leq s \leq t}|M_s^{n+1, m}|.\] 
The second assertion is satisfied for $m=0$. Suppose $(ii)$ holds for all $k\leq m$ and let us prove it for $k=m+1$. For this task we use again Proposition \ref{prop1:Nualart}. Using (\ref{eq3:appro}), we first need to prove that $\max_{0\leq u\leq s}(-a^{n}_u -\alpha M^{n+1, m}_u) \in \mathbb{D}^{1, 2}$. Indeed, we have   $a^{n}_u -\alpha M^{n+1, m}_u \in \mathbb{D}^{1, 2}$ and 
\[ \E[\sup_{0\leq u \leq s }(-a^{n}_u -\alpha M^{n+1, m}_u)^2]\leq 2 \E[\sup_{0\leq u \leq s }(a^{n}_u )^2] + 2 \alpha^2 \E[\sup_{0\leq u \leq s }( M^{n+1, m}_u )^2] <\infty.  \]
Furthermore, using $(\ref{eq:max an})$ and the induction hypothesis for $m$, we get
\begin{align*}
\E[\sup_{0\leq u \leq s }||D(-a^{n}_u -\alpha M^{n+1, m}_u)||_H^2] \leq 2 \E[\sup_{0\leq u \leq s }||Da^{n}_u||_H^2] + 2 \E[\sup_{0\leq u \leq s }||D(M^{n+1, m}_u)||^2_H] <\infty.
\end{align*}

Therefore, by Proposition \ref{prop1:Nualart}  and the point $(i)$ of Lemma \ref{lem1} we deduce that
\begin{align} \label{aux1:lem}
a^n_s+ \sup_{0\leq u\leq s}(-a^{n}_u -\alpha M^{n+1, m}_u) \in \mathbb{D}^{1, 2}.
\end{align}
Moreover,  
\begin{align} \label{aux2:lem}
&\E\left[\sup_{0 \leq s \leq t}\left(a^n_s +  \frac{\beta}{\beta -1} \max_{0\leq u\leq s }(-a^n_u-\alpha M^{n+1, m}_u)\right)^2\right]\nonumber\\
& \leq 2\left(1 +\frac{|\beta|^2}{(1-\beta)^2}\right) \E\left[\max_{0\leq u\leq s }(a^n_u)^2\right] + \frac{2 \alpha^2 \beta^2}{(1-\beta)^2} \E\left[\sup_{0\leq u \leq t}(M ^{n+1, m}_u)^2\right]<\infty
\end{align}
and  we have\footnote{One can easily check the assumptions of Proposition \ref{prop1:Nualart}.}
\[||  D(\max_{0\leq u\leq s}(-a^{n}_u -\alpha M^{n+1, m}_u)) ||_H^2 \leq  \max_{0\leq u\leq s}||D(-a^{n}_u -\alpha M^{n+1, m}_u)) ||_H^2\, \, \, \mbox{a.e.} \]
so, a.e.
\begin{align*}
 ||D(-a^{n}_s +  \frac{\beta}{\beta -1} \max_{0\leq u\leq s }(-a^n_u-\alpha M^{n+1, m}_u))||^2_H &\leq (2+ \frac{4 \beta^2}{(1-\beta)^2}) \max_{0\leq s \leq t} || Da^{n}_s||^2_H\\
& + \frac{4  \alpha^2\beta^2}{(1-\beta)^2} \max_{0\leq s\leq t}||DM^{n+1, m} ||^2_H 
\end{align*}
and  then 
\begin{align}\label{aux3:lem}
 &\E\left[\max_{0\leq s \leq t}||D(-a^{n}_s +  \frac{\beta}{\beta -1} \max_{0\leq u\leq s }(-a^n_u-\alpha M^{n+1, m}_u))||^2_H\right] \nonumber\\
 & \leq \left(2+ \frac{4 \beta^2}{(1-\beta)^2}\right) \E[\max_{0\leq s \leq t} || Da^{n}_s||^2_H] + \frac{4  \alpha^2\beta^2}{(1-\beta)^2} \E[ \max_{0\leq s\leq t}||DM_s^{n+1, m} ||^2_H ],
\end{align}
which is finite by  (\ref{eq:max an}) in Lemma \ref{lem1} and induction hypothesis on $m$.

\vspace{0.1cm}

\noindent Therefore (\ref{aux1:lem}), (\ref{aux2:lem}) and  (\ref{aux3:lem}) together with Proposition \ref{prop1:Nualart} imply
\[M_t^{n+1, m+1} =\frac{1}{1-\alpha}\max_{0\leq s \leq t}\left(a^n_s + \frac{\beta}{\beta -1} \max_{0\leq u\leq s }(-a^n_u-\alpha M^{n+1, m}_u)\right) \in \mathbb{D}^{1, 2}.\]
This finishes the proof of the first claim of $(ii)$ in Lemma \ref{lem2}. The second claim is easily deduced since, again by Proposition \ref{prop1:Nualart},  one has a.e.
\begin{align}\label{ineq:1}
&||DM_s^{n+1, m+1}||^2_H \leq \frac{1}{(1-\alpha)^2}\max_{0\leq u\leq s}|| D (a^n_u +  \frac{\beta}{\beta -1} \max_{0\leq v\leq u }(-a^n_v-\alpha  M^{n+1, m}_v))||^2_H \nonumber\\
&\leq \frac{2}{(1-\alpha)^2}\left(1+ \frac{\beta^2}{(1-\beta)^2}\right) \max_{0\leq u\leq s}|| Da^n_u||^2_H + \frac{2\alpha^2\beta^2}{(1-\alpha)^2(1-\beta)^2}\max_{0\leq u\leq s }||DM_s^{n+1, m}||^2_H.
\end{align}
Next, let  us prove the last point $(iii)$ of Lemma \ref{lem2}. A more careful analysis of (\ref{ineq:1}), using the elementary inequality $(a+b)^2\leq (1+\varepsilon^{-1})a^2 + (1+\varepsilon)b^2$ valid for all real $a$, $b$ and $\varepsilon>0$, leads to 
\begin{align}\label{ineq:2}
&||DM_s^{n+1, m+1}||^2_H \nonumber\\
&\leq  \frac{1+\varepsilon^{-1}}{(1-\alpha)^2}\max_{0\leq u\leq s}|| D(a^n_u)||^2_H + \frac{(1+\varepsilon)\beta^2}{(1-\beta)^2(1-\alpha)^2}\max_{0\leq u\leq s}|| D(-a^n_u-\alpha  M^{n+1, m}_u))||^2_H \nonumber\\
&\leq \frac{1+\varepsilon^{-1}}{(1-\alpha)^2}\left(1+ \frac{(1+\varepsilon)\beta^2}{(1-\beta)^2}\right)\max_{0\leq u\leq s}|| D(a^n_u)||^2_H + \frac{(1+\varepsilon)^2 \alpha^2\beta^2}{(1-\alpha)^2(1-\beta)^2}\max_{0\leq u\leq s }||DM_s^{n+1, m}||^2_H\, \mbox{a.e.}
\end{align}
and then 
\begin{align}\label{ineq:3}
& \E[\max_{0\leq s\leq t }||DM_s^{n+1, m}||^2_H] \leq \frac{1+\varepsilon^{-1}}{(1-\alpha)^2}\left(1+ \frac{(1+\varepsilon)\beta^2}{(1-\beta)^2}\right)\E[\max_{0\leq s\leq t}|| D(a^n_s)||^2_H] \nonumber\\
&+ \frac{(1+\varepsilon)^2 \alpha^2\beta^2}{(1-\alpha)^2(1-\beta)^2}\E[\max_{0\leq s\leq t }||DM_s^{n+1, m}||^2_H].
\end{align}
By choosing $\varepsilon >0 $ so that $\frac{(1+\varepsilon)^2 \alpha^2\beta^2}{(1-\alpha)^2(1-\beta)^2}<1$  and iterating (\ref{ineq:3}), we deduce  the claim. This finishes the proof of Lemma \ref{lem2}.\\
\noindent {\bf Continuation of the proof of Lemma \ref{lem1}:}  Thanks to (\ref{conv:L2}) and Lemma \ref{lem2} one can use  Proposition \ref{appro:D12} to deduce that $M^n_t  \in \mathbb{D}^{1, 2}$. Next, recall first from \cite{Belfadli2009, CD1999} that
\begin{equation}\label{eq4}
M_t^{n+1} =\frac{1}{1-\alpha}\max_{0\leq s \leq t}\left\{a^n_s + \frac{\beta}{\beta -1} \max_{0\leq u\leq s }(-a^n_u-\alpha M^{n+1}_u)\right\}.  
\end{equation}

Then, arguing as above and choosing $\varepsilon>0$ so that   $\frac{(1+\varepsilon)^2 \alpha^2\beta^2}{(1-\alpha)^2(1-\beta)^2}<1$  we obtain
\begin{align}\label{ineq:4}
& \left(1-\frac{(1+\varepsilon)^2 \alpha^2\beta^2}{(1-\alpha)^2(1-\beta)^2}\right)\E[\max_{0\leq s\leq t }||DM_s^{n+1}||^2_H]\nonumber\\
& \leq \frac{1+\varepsilon^{-1}}{(1-\alpha)^2}\left(1+ \frac{(1+\varepsilon)\beta^2}{(1-\beta)^2}\right)\E[\max_{0\leq s\leq t}|| D(a^n_s)||^2_H]. 
\end{align}
This proves the second claim $(ii)$ of Lemma \ref{lem1}.\\  Concerning  $I^n_t$ we recall from \cite{Belfadli2009} that
\begin{align}
&(\be - 1)I_t^{n+1} =\disp{\max_{0\leq s\leq
t}{\left(-a^n_s- \al
M_s^{n+1}\right)}}\label{In:1}\\
&(1-\al)M_t^{n+1} =\disp{\max_{0\leq s\leq t}{\left(a^n_s+ \be I^{n+1}_s\right)}}\label{In:2}.
\end{align}
Combining these two equalities, we get:
\begin{eqnarray}\label{e6}
I_t^{n+1}&=&\frac{1}{\beta- 1}\max_{0\leq s\leq
t}\left\{-a^{n}_s -\frac{\alpha}{1- \alpha}\max_{0\leq u\leq s}(a^n_s + \beta I_u^{n+1})\right\}.
\end{eqnarray}
On one hand, using the fact that $M^{n+1}_s \in \mathbb{D}^{1, 2}$,  (\ref{In:1}) and Proposition \ref{prop1:Nualart}, we obtain without difficulty that $I^{n+1}_t$ belongs to $\mathbb{D}^{1, 2}$. On the other hand, arguing as in (\ref{ineq:4})  by using $(\ref{e6})$ we obtain, for $\varepsilon>0$ such that  $\frac{(1+\varepsilon)^2 \alpha^2\beta^2}{(1-\alpha)^2(1-\beta)^2}<1$ ,
\begin{align}\label{ineq:5}
& \left(1-\frac{(1+\varepsilon)^2 \alpha^2\beta^2}{(1-\alpha)^2(1-\beta)^2}\right)\E[\max_{0\leq s\leq t }||DI_s^{n+1}||^2_H]\nonumber\\
& \leq \frac{1+\varepsilon^{-1}}{(1-\beta)^2}\left(1+ \frac{(1+\varepsilon)\alpha^2}{(1-\alpha)^2}\right)\E[\max_{0\leq s\leq t}|| D(a^n_s)||^2_H]. 
\end{align} 
This proves the third claim of Lemma 3.5 and finishes the proof of Lemma  \ref{lem2}.
\subsection{Proof of Theorem $\ref{thm1}$}
\noindent {\bf Proof of Theorem  \ref{thm1}}: We will use Proposition $\ref{appro:D12}$. Taking into account Proposition \ref{proposition1} and the convergence $(\ref{convergence1})$, we just need to show that 
\[
\sup_{n} \E\left[||DX^{n}_t ||_{H}^2\right]<\infty.\]
For every  $r \in[0, t]$ and integer $n$,
$$
\begin{aligned}
 D_r X_t^{n+1}=\sigma\left(X_r^n\right)+\int_r^t D_r\left(\sigma\left(X_s^n\right)\right) d W_s+\int_r^t D_r\left(b\left(X_s^n\right)\right) ds+ \alpha D_rM_t^{n+1} +  \beta D_r I_t^{n+1}. 
\end{aligned}
$$
Using the elementary inequality  $|a+b+c+d+e|^2 \leqslant 5\left(|a|^2+|b|^2+|c|^2+|d|^2+|e|^2\right)$ valid for any real numbers $a, b, c, d$ and $e$, we get :
$$
\begin{aligned}
& \E[||DX^{n+1}_t||^2_{H}]= \E\left[ \int_{0}^t |D_rX^{n+1}_t|^2dr\right]\\
&\leq  5 \left\{  \E \left[ \int_{0}^t |\sigma(X^n_r)|^2 dr \right] + \E \left[ \int_{0}^t \left| \int_{r}^t D_r(\sigma(X^n_s))dW_s\right|^2 dr \right]\right.
\\
& \left. + \E \left[ \int_{0}^t \left| \int_{r}^t D_r(b(X^n_s))ds \right|^2 dr \right] + \alpha^2 \E\left[\left|\int_{r}^{t} D_rM^{n+1}_t dr\right|^2\right] +  \beta^2 \E\left[\left|\int_{r}^{t} D_rI^{n+1}_t dr\right|^2\right] \right\}\\
&=: I^{n, 1}_t + I^{n, 2}_t + I^{n, 3}_t + I^{n, 4}_t + I^{n, 5}_t. 
\end{aligned}
$$
Let us estimate each term of this inequality. For $I^{n, 1}_t$, by Lipschitz condition on $\sigma$, we have
\begin{align}
I^{n, 1}_t \leq c \int_{0}^t \E(1+ |X^{n}_r|^2) dr \leq c  \sup_{n} \int_{0}^t \E(1+ |X^{n}_r|^2) dr=: A,
\end{align}
which is finite since $(X^n)$ converges to $X$ uniformly with respect to time parameter, see for instance \cite{Belfadli2009}.  \\
\noindent For  the terms $I^{n, 2}_t$ and $I^{n, 3}_t$, using (3), (4), It\^{o}'s isometry and Cauchy-Schwarz inequality we obtain
\begin{align}
I^{n, 2}_t + I^{n, 3}_t \leq 5 K^2 (1 +t) \int_{0}^t \int_{r}^{t}\E[ |D_rX^{n}_s|^2]ds dr.
\end{align}
For the term $I^{n, 4}_t$ ($I^{n, 5}_t$ can be treated in a similar way), using Cauchy-Schwarz inequality we have
\[I^{n, 4}_t \leq \alpha^2 t \E[||DM^{n+1}_t||^2_{H}].\]
Then, using the identity $(\ref{eq4})$ we have a.e.
\begin{align*}
% \nonumber to remove numbering (before each equation)
&||DM^{n+1}_t||^2_{H} = \frac{1}{(1-\alpha)^2}||D\left( \max_{0\leq s \leq t} \left(a^n_s + \frac{\beta}{\beta -1} \max_{0\leq u \leq s}(-a^{n}_u -\alpha M^{n+1}_u)\right)\right)||_{H}^2\\
& \leq   \frac{1}{(1-\alpha)^2}\max_{0\leq s \leq t}  ||D \left(a^n_s + \frac{\beta}{\beta -1} \max_{0\leq u \leq s}(-a^{n}_u -\alpha M^{n+1}_u)\right)||_{H}^2 \\
& \leq \frac{1+\varepsilon^{-1}}{(1-\alpha)^2} \max_{0 \leq s \leq t} || Da^{n}_s||^2_{H}  +  \frac{(1 +\varepsilon)\beta^2}{(1-\alpha)^2(1-\beta)^2} \max_{0\leq s \leq t }||D\max_{0\leq u \leq s }( -a^{n}_u -\alpha M^{n+1}_u) ||_{H}^2\\
& \leq \frac{1+\varepsilon^{-1}}{(1-\alpha)^2} \max_{0 \leq s \leq t} || Da^{n}_s||^2_{H}  +   \frac{(1 +\varepsilon)(1 +\varepsilon^{-1})\beta^2}{(1-\alpha)^2(1-\beta)^2} \max_{0\leq s \leq t }|| Da^{n}_s||^2_{H} \\
& + \frac{(1+ \varepsilon)^2 \alpha^2\beta^2}{(1-\alpha)^2(1-\beta)^2} \max_{0\leq s \leq t}|| DM^{n+1}_s||^2_{H}.
\end{align*}
So that a.e.
\begin{align*}
\left ( 1 -\frac{(1+ \varepsilon)^2 \alpha^2\beta^2}{(1-\alpha)^2(1-\beta)^2}\right)  \max_{0 \leq s \leq t} || DM^{n+1}_s||^2_{H} \leq \left( \frac{1+\varepsilon^{-1}}{(1-\alpha)^2} + \frac{(1+ \varepsilon)(1+ \varepsilon^{-1})\beta^2}{(1-\alpha)^2(1-\beta)^2}\right) \max_{0 \leq s \leq t} || Da^{n}_s||^2_{H},  
\end{align*}
where $\varepsilon$ is chosen such that $\frac{(1+ \varepsilon)^2 \alpha^2\beta^2}{(1-\alpha)^2(1-\beta)^2} <1$. Therefore
\begin{align*}
\sup_{0\leq s \leq t} \E\left(  || DM^{n+1}_s||^2_{H}\right) \leq \E (\sup_{0\leq s \leq t}    || DM^{n+1}_s||^2_{H} ) \leq c \E(\max_{0 \leq s \leq t}  || Da^{n}_s||^2_{H}).
\end{align*}
Similarly, for $I^{n+1}_t$, using  $(\ref{e6})$ we obtain a.e.
\begin{align*}
\left ( 1 -\frac{(1+ \varepsilon)^2 \alpha^2\beta^2}{(1-\alpha)^2(1-\beta)^2}\right)  \max_{0 \leq s \leq t} || DI^{n+1}_s||^2_{H} \leq \left( \frac{1+\varepsilon^{-1}}{(1-\beta)^2} + \frac{(1+ \varepsilon)(1+ \varepsilon^{-1})\alpha^2}{(1-\alpha)^2(1-\beta)^2}\right) \max_{0 \leq s \leq t} || Da^{n}_s||^2_{H},  
\end{align*}
and then 
\begin{align*}
\sup_{0\leq s \leq t} \E\left[ || DI^{n+1}_s||^2_{H}\right] \leq \E [\sup_{0\leq s \leq t}    || DI^{n+1}_s||^2_{H} ] \leq c \E[\max_{0 \leq s \leq t}  || Da^{n}_s||^2_{H}].
\end{align*}
Implying in turns that 
\begin{align}
I^{n, 4}_t +I^{n, 5}_t \leq c \E[\max_{0 \leq s \leq t}  || Da^{n}_s||^2_{H}].
\end{align}
 On the other hand,  using the estimation $(\ref{max:e2})$, we have 

\[ \E[\max_{0 \leq s \leq t}  || Da^{n}_s||^2_{H}] \leq  c\left(A +\int_{0}^t \int_{r}^t \E[|D_rX^{n}_{u}|^2]dudr\right).\]
Therefore
\[I^{n, 4}_t +I^{n, 5}_t \leq c \left(A +\int_{0}^t \int_{r}^t \E[|D_rX^{n}_{u}|^2]dudr\right),\]
and then 
\[ \E(||DX^{n+1}_t||^2_{H}) \leq c\left(1 + \int_{0}^t  \E[||DX^{n}_s||^2_{H}]ds\right).\] 
Iterating this inequality, we obtain $\sup_{n} \E\left[||DX^{n}_t ||_{H}^2\right]<\infty.$ Thus by Proposition $\ref{appro:D12}$ $X_t \in \mathbb{D}^{1, 2}$ and $DX^{n}_t$  converges weakly to $DX_t$ in $L^{2}(\Omega; H)$ and  finishing the proof of Theorem 1.

A further generalization of Theorem \ref{thm1} may be stated as follows.
\begin{prop} \label{generalization of thm1} Let $\alpha <1$ and $\beta<1$ be such that $|\rho| <1$. 
Let $(X_t)_{t\geq 0}$ be the unique solution of Eq. (\ref{e2}). Then,  $X_t $ belong to $\mathbb{D}^{k, 2}$ for any $t\geq 0$ and integer $k\geq 1$.  
\end{prop}
\noindent {\bf Proof }:  Using the Picard's iteration scheme $(X^n)_n$ introduced in (\ref{eq2}), the proof is an immediate consequence of Proposition \ref{appro:D12}. Indeed, proceeding as in the proof of Theorem \ref{thm1} by using  the sequence  $(M^{n, m})_m$ and  an  induction argument on $n$ and $k$, one can check that 

\[X^n \in \mathbb{D}^{k, 2}\,\, \mbox{and}\, 
\sup_{n} \E\left[||D^{i}X^{n}_t ||_{{H^{\otimes i}}}^2\right]<\infty, \, \forall i= 1, \cdots, k.\]

\section{Absolute continuity of the laws of doubly perturbed diffusion processes}
In this section, assuming Lipschitz condition on $b$ and $\sigma$, and a non degeneracy condition on $\sigma$, we prove the absolute continuity of the laws of $X_t$, for all $t>0$. 
\begin{thm} \label{thm2} Let $\alpha <1$ and $\beta<1$ such that $|\rho| <1$. 
Let $(X_t)_{t\geq 0}$ be the unique solution of Eq. (\ref{e2}).  Assume that $b$ and $\sigma$ are Lipschitz continuous  and $|\sigma(x)| >0$, for all $x\in \mathbb{R}$. Then, for all $t>0$, the law of $X_t$ is absolutely continuous  with respect to the Lebesgue measure.
\end{thm}

\noindent {\bf Proof of Theorem \ref{thm2}}: In virtue  of  Proposition \ref{absl}, it is enough to show that $||DX_t ||_H>0$ a.s.
\noindent For $r\leq t$, we have 
 \begin{align}\label{derivative eq}
& D_rX_t= \sigma(X_r)+ \int_{r}^{t}D_{r}(\sigma(X_s))dW_s + \int_{r}^{t}D_{r}(b(X_s))ds\nonumber\\
& + \alpha D_r(\max_{0\leq s \leq t}X_s) + \beta D_r(\min_{0\leq s \leq t}X_s). 
\end{align}
Using the elementary inequality $(a+b)^2\geq \frac{1}{2}a^2 -b^2$, we have 
 \begin{align}
 &(D_rX_t)^2\nonumber\\
 &\geq \frac{1}{2} \sigma(X_r)^2 -\left( \int_{r}^{t}D_{r}(\sigma(X_s))dW_s + \int_{r}^{t}D_{r}(b(X_s))ds + \alpha D_r(\max_{0\leq s \leq t}X_s) + \beta D_r(\min_{0\leq s \leq t}X_s)\right)^2\nonumber \\
& \geq  \frac{1}{2} \sigma(X_r)^2 -4 \left(\int_{r}^{t}D_{r}(\sigma(X_s))dW_s\right)^2 -4 \left( \int_{r}^{t}D_{r}(b(X_s))ds\right)^2 \nonumber\\
& -4 \alpha^2 (  D_r(\max_{0\leq s \leq t}X_s))^2-4 \beta^2 (  D_r(\min_{0\leq s \leq t}X_s))^2.
\end{align}
The continuity of $r\mapsto \sigma(X_r)$ gives
\begin{align}\label{first term}
\lim_{\varepsilon \rightarrow 0}\frac{1}{\varepsilon} \int_{t-\varepsilon}^t \sigma (X_r)^2 dr = \sigma(X_t)^2.
\end{align}
For the second term, we have by Cauchy-Schwarz inequality and Fubini's Theorem,  for all $0<\varepsilon< t/2$,  
\begin{align}\label{second term}
&\E \left[\int_{t-\varepsilon}^t (\int_{r}^t D_r(b(X_s))ds)^2 dr\right] \leq \int_{t-\varepsilon}^t (t-r) \int_r^t \E[D_r(b(X_s))^2]dsdr \nonumber\\
& \leq K^2\varepsilon \int_{t-\varepsilon}^t \int_{t-\varepsilon}^s \E[(D_rX_s)^2] drds \leq K^2\varepsilon \int_{t-\varepsilon}^t \int_{s-\varepsilon}^s \E[(D_rX_s)^2] drds.
\end{align}
For the third term $\E[ \int_{t-\varepsilon}^t (\int_{r}^t D_r(\sigma(X_s))dW_s)^2 dr]$, we have
\begin{align}\label{ineq:third term} 
& \E \left[ \int_{t-\varepsilon}^t (\int_{r}^t D_r(\sigma(X_s))dW_s)^2 dr\right] \leq K^2 \int_{t-\varepsilon}^t \int_r^t \E[(D_rX_s)^2]ds dr \nonumber \\
& \leq K^2 \int_{t-\varepsilon}^t \int_{s-\varepsilon}^s \E[(D_rX_s)^2]dr ds.
\end{align}
Let us  first show the following claim
\begin{claim} \label{claim2} We have, for all $ s \in [t-\varepsilon, t]$, 
\begin{align} 
\int_{s-\varepsilon}^s \E[(D_rX_s)^2] dr \leq c \varepsilon.
\end{align}
\end{claim}

\noindent {\bf Proof of Claim \ref{claim2}}:  We  use the sequence of approximations $(X^n)$ in (\ref{eq2}).

\noindent For $r\leq t$, we have 
 \begin{align*}
& D_rX^{n+1}_t= \sigma(X^{n}_r)+ \int_{r}^{t}D_{r}(\sigma(X^{n}_s))dW_s + \int_{r}^{t}D_{r}(b(X_s^n))ds\\
& + \alpha D_r(M^{n+1}_s) + \beta D_r(I^{n+1}_s). 
\end{align*}
Then 
 \begin{align}\label{iteration1}
& \varphi_{n+1}(s):= \int_{s-\varepsilon}^s \E[(D_rX^{n+1}_s)^2] dr \leq 5 \E \int_{s-\varepsilon}^s \sigma( X^n_r)^2 dr + 5 \E[ \int_{s-\varepsilon}^s dr (\int_r^s D_r(\sigma(X^n_u)) dW_u)^2]\nonumber\\
& +5\E[ \int_{s-\varepsilon}^s dr (\int_r^s D_r(b(X^n_u)) du)^2] + 5 \alpha^2 \E[\int_{s-\varepsilon}^s (D_rM^{n+1}_s)^2 dr] + 5 \beta^2 \E[\int_{s-\varepsilon}^s (D_rI^{n+1}_s)^2 dr]\nonumber\\
&=: \sum_{i=1}^5 J^{n, i}_{\varepsilon}(s).
\end{align}
We have 
\begin{align}\label{first terms}
\sum_{i=1}^3 J^{n, i}_{\varepsilon}(s) \leq 5 \int_{s-\varepsilon}^s \E[\sigma(X^n_r)^2]dr + 10 K^2 \E\left[\int_{s-\varepsilon}^s\int_{r}^s (D_rX^n_u)^2 dudr\right].
\end{align}
For $J^{n, 4}_{\varepsilon}(s)$ and $J^{n, 5}_{\varepsilon}(s)$, notice first that, by a slight modification of the argument used in the first section, we have 
for $\gamma >0$ such that $\frac{(1+\gamma)^2 \alpha^2 \beta^2}{(1-\alpha)^2(1-\beta)^2}<1$, 
\begin{align}\label{modi1}
& \left(1-\frac{(1+\gamma)^2 \alpha^2\beta^2}{(1-\alpha)^2(1-\beta)^2}\right)\E\left[\int_{s-\varepsilon}^s(D_rM^{n+1}_s)^2dr\right]\nonumber\\
&\leq  \frac{1+\gamma^{-1}}{(1-\alpha)^2}\left(1+ \frac{(1+\gamma)\beta^2}{(1-\beta)^2}\right)\E\left[\max_{0\leq u\leq s}\left( \int_{s-\varepsilon}^s |D_ra^n_u|^2dr\right)\right]. 
\end{align}
 and 
 \begin{align}\label{modi2}
& \left(1-\frac{(1+\gamma)^2 \alpha^2\beta^2}{(1-\alpha)^2(1-\beta)^2}\right)\E\left[\int_{s-\varepsilon}^s(D_rI^{n+1}_s)^2dr\right]\nonumber\\
&\leq  \frac{1+\gamma^{-1}}{(1-\beta)^2}\left(1+ \frac{(1+\gamma)\alpha^2}{(1-\alpha)^2}\right)\E\left[\max_{0\leq u\leq s}\left( \int_{s-\varepsilon}^s |D_ra^n_u|^2dr\right)\right]. 
\end{align}
On the other hand, using (\ref{eq:a}) it is easily seen that 
\begin{align} \label{mod:max:Dan}
&\E\left[\max_{0\leq u\leq s}\left( \int_{s-\varepsilon}^s |D_ra^n_u|^2dr\right)\right] \leq  3 \int_{s-\varepsilon}^s  \E[\sigma(X^n_r)^2] dr +(3K^2+cK^2)\int_{s-\varepsilon}^s \int_{r}^s \E[|D_rX^{n}_v|^2]dv dr.
\end{align}
Then $(\ref{modi1})$, $(\ref{modi2})$ and  $(\ref{mod:max:Dan})$, together yield
\begin{align}\label{firth-fifth term}
J^{n, 4}_{\varepsilon}(s) +J^{n, 5}_{\varepsilon}(s) \leq c\int_{s-\varepsilon}^s \E[\sigma(X^n_r)^2] dr +c^{\prime} \int_{s-\varepsilon}^s \int_{r}^s \E[|D_rX^{n}_v|^2]dv dr.
\end{align}
 Combining  (\ref{first terms}),  (\ref{firth-fifth term}) and  (\ref{iteration1}) give
 \begin{align*}
 &\varphi_{n+1}(s) \leq c  \int_{s-\varepsilon}^s \E[\sigma(X^n_r)^2] dr  + c^{\prime}  \int_{s-\varepsilon}^s  dr \int_{r}^sdv  \E[|D_rX^{n}_v|^2]\\
 & \leq c  \int_{s-\varepsilon}^s \E[\sigma(X^n_r)^2] dr  + c^{\prime} \int_{s-\varepsilon}^s  dv \int_{v-\varepsilon}^v dr \E[|D_rX^{n}_v|^2]. 
 \end{align*}
 That is 
 \[ \varphi_{n+1}(s) \leq c  \int_{s-\varepsilon}^s \E[\sigma(X^n_r)^2] dr  + c^{\prime} \int_{s-\varepsilon}^s  \varphi_{n}(v) dv. 
 \]
 But 
 \[ \int_{s-\varepsilon}^s \E[\sigma(X^n_r)^2] dr \leq K \int_{s-\varepsilon}^s (1+ \E[(X^{n}_{r})^2])dr \leq (1 + \sup_n \E[\sup_{0\leq r \leq s} (X^{n}_{r})^2] )K \varepsilon= c\varepsilon.
  \]
  Therefore
  \begin{align}
 \varphi_{n+1}(s) \leq  c\varepsilon  +  c^{\prime} \int_{s-\varepsilon}^s  \varphi_{n}(v) dv. 
 \end{align}
 Iterating this inequality we get 
 \begin{align*}
 &  \varphi_{n}(s) \leq c \varepsilon \sum_{j=0}^{n-1} (c^{\prime} \varepsilon)^j  + c^{\prime \prime} \varepsilon^{n}\leq c^{\prime \prime \prime} \varepsilon \sum_{j=0}^{n-1}(c^{\prime \prime \prime} \varepsilon)^j \leq c^{\prime \prime \prime} \varepsilon \frac{1}{1-c^{ \prime \prime \prime} \varepsilon} \leq 2 c^{\prime \prime \prime} \varepsilon,
 \end{align*}
 by choosing $\varepsilon < \frac{1}{2 c^{\prime \prime \prime}}$. 
 Finally, weak convergence of $DX^{n}$ to $DX$ imply 
 \[\E\left[\int_{s-\varepsilon}^s (D_rX_s)^2 dr\right]\leq   \liminf_{n\rightarrow +\infty} \varphi_n(s) \leq 2 c^{\prime \prime \prime} \varepsilon.\]
 This finishes the proof of the claim.\\
Now, applying Claim \ref{claim2} to (\ref{second term}) and (\ref{ineq:third term}), we  conclude  that
 \begin{align}\label{eq:1 and 2}
&\lim_{\varepsilon \rightarrow 0}\frac{1}{\varepsilon} \left\{\E\left[ \int_{t-\varepsilon}^t \left(\int_{r}^t D_r(\sigma(X_s))dW_s\right)^2 dr\right] +  \E\left[ \int_{t-\varepsilon}^t \left(\int_{r}^t D_r(b(X_s))ds\right)^2 dr \right] \right\}=0.
\end{align}
As a by product, we derive that there exists a sequence $(\varepsilon_n)_n \downarrow 0$ and  a negligible set $\mathcal{N}$ such that  on $\mathcal{N}^{c}$
  \begin{align}\label{eq:1 and 2 a.s}
&\lim_{n  \rightarrow +\infty}\frac{1}{\varepsilon_n}  \int_{t-\varepsilon_n}^t \left[\left(\int_{r}^t D_r(\sigma(X_s))dW_s\right)^2 + \left(\int_{r}^t D_r(b(X_s))ds\right)^2\right] dr =0.
\end{align}
Set, for any integer  $n$,
\begin{align*}
& A_n:=\{\omega \in \Omega: \max_{0\leq s \leq t}X_s(\omega)=\max_{0\leq s \leq t-\varepsilon_n}X_s(\omega) \},\, \,  A:=\{\omega \in \Omega: \max_{0\leq s \leq t}X_s(\omega)=X_t(\omega) \}
\end{align*}
 and similarly
 \begin{align*}
& B_n:=\{\omega \in \Omega: \min_{0\leq s \leq t}X_s(\omega)=\min_{0\leq s \leq t-\varepsilon_n}X_s(\omega) \},  B:=\{\omega \in \Omega: \min_{0\leq s \leq t}X_s(\omega)=X_t(\omega) \}.
\end{align*}
Clearly  $\Omega= \bigcup_{n, m} (A_n\cap B_m) \cup \bigcup_{n} (A_n\cap B) \cup \bigcup_{m} (B_m\cap A) \cup (A\cap B).$\\

$\bullet$ If $\omega \in \mathcal{N}^{c} \cap A_n \cap B_m$ for some integers  $n$ and $m$, then for any $j> n+m$, we have
\[\alpha^2 \int_{t-\varepsilon_j}^t (D_r(\max_{0\leq s\leq t-\varepsilon_n}X_s(\omega)))^2dr=0\]
and 
\[\beta^2 \int_{t-\varepsilon_j}^t (D_r(\min_{0\leq s\leq t-\varepsilon_m}X_s(\omega)))^2dr=0.\]
Hence,  by the local property of the Malliavin derivative operator (see for instance page $49$ of \cite{Nualart2006}), we obtain 
\begin{align}
&\lim_{j  \rightarrow +\infty}\frac{1}{\varepsilon_j}  \int_{t-\varepsilon_j}^t (D_r(\max_{0\leq s\leq t} X_s(\omega)))^2dr = \lim_{j  \rightarrow +\infty}\frac{1}{\varepsilon_j}  \int_{t-\varepsilon_j}^t (D_r(\max_{0\leq s\leq t-\varepsilon_n}X_s(\omega)))^2dr=0,
\end{align}
 and
\begin{align*}
&\lim_{j  \rightarrow +\infty}\frac{1}{\varepsilon_j}  \int_{t-\varepsilon_j}^t (D_r(\min_{0\leq s\leq t} X_s(\omega)))^2dr = \lim_{j  \rightarrow +\infty}\frac{1}{\varepsilon_j}  \int_{t-\varepsilon_j}^t (D_r(\min_{0\leq s\leq t-\varepsilon_m}X_s(\omega)))^2dr=0,
\end{align*}
 on $\mathcal{N}^c\cap A_n\cap B_m$. Therefore, for all $\omega \in \bigcup_{n, m} (A_n\cap B_m)$: 
\begin{align}
&\lim_{j  \rightarrow +\infty}\frac{1}{\varepsilon_j}  \int_{t-\varepsilon_j}^t (D_rX_t)^2 dr \geq \frac{1}{2} \sigma(X_t)^2>0 \, \, \, a.s.
\end{align}

$\bullet$  For $\omega \in \mathcal{N}^{c} \cap A\cap B_m$ for some integer  $m$,  we have 
\[(1-\alpha) D_rX_t =\sigma(X_r) + \int_{r}^{t}D_{r}(\sigma(X_s))dW_s + \int_{r}^{t}D_{r}(b(X_s))ds  + \beta D_r(\min_{0\leq s \leq t -\varepsilon_m}X_s)\]

then for all $j> m$, 
\begin{align*}
(1-\alpha)^2 \frac{1}{\varepsilon_j}  \int_{t-\varepsilon_j}^t (D_rX_t)^2 dr& \geq \frac{1}{2} \sigma(X_t)^2  - 2 \int_{t -\varepsilon_j}^{t}\left(\int_{r}^t D_r(\sigma(X_s))dW_s\right)^2 dr\\
& - 2 \int_{t -\varepsilon_j}^{t}\left(\int_{r}^t D_r(b(X_s))ds\right)^2 dr-2 \beta^2   \int_{t-\varepsilon_j}^t (D_r(\min_{0\leq s \leq t -\varepsilon_m}X_s))^2 dr. 
\end{align*}
But 
\[
 \lim_{j  \rightarrow +\infty} \frac{1}{\varepsilon_j}\int_{t-\varepsilon_j}^t (D_r(\min_{0\leq s \leq t -\varepsilon_m}X_s))^2 dr=0\]
hence using (\ref{eq:1 and 2 a.s}) we deduce
\begin{align}
&\lim_{j  \rightarrow +\infty}\frac{(1-\alpha)^2}{\varepsilon_j}  \int_{t-\varepsilon_j}^t (D_rX_t)^2 dr \geq \frac{1}{2} \sigma(X_t)^2 > 0, \, \, \, \, \mbox{a.s. on }\,\,   \bigcup_{m} (A\cap B_m).
\end{align}

$\bullet$  For $\omega \in  \bigcup_{m} (A\cap B_m)\cap  \mathcal{N}^{c}$, similar arguments as before lead to 
\begin{align}
&\lim_{j  \rightarrow +\infty}\frac{(1-\beta)^2}{\varepsilon_j}  \int_{t-\varepsilon_j}^t (D_rX_t)^2 dr \geq \frac{1}{2} \sigma(X_t)^2>0.
\end{align}

$\bullet$  For $\omega \in  A\cap B\cap \mathcal{N}^{c}$,  (\ref{derivative eq}) gives
\[(1-\alpha -\beta)^2 (D_rX_t)^2 \geq \frac{1}{2} \sigma(X_r)^2 -(\int_{r}^{t}D_{r}(\sigma(X_s))dW_s + \int_{r}^{t}D_{r}(b(X_s))ds)^2.\]
So on $A\cap B\cap \mathcal{N}^{c}$, 
\begin{align}
&\lim_{j  \rightarrow +\infty}\frac{1}{\varepsilon_j}  \int_{t-\varepsilon_j}^t (D_rX_t)^2 dr \geq \frac{1}{2(1-\alpha -\beta)^2} \sigma(X_t)^2>0 \, \, \, a.s.
\end{align}
Therefore 
\[|| DX_t||^2_H =\int_{0}^t (D_rX_t)^2 dr >0 \,\, \mbox{a.s.}\]
Thus finishing the proof of Theorem \ref{thm2}.
\section{Smoothness of the density}  
 
 The purpose of this section is to study the smoothness of the density for   each $X_t$, where  $X$ denotes the unique solution  of Eq. (\ref{e2}).

As  previously  we assume that $\alpha <1$ and $\beta<1$ are such that $|\rho| <1$. We also assume that $b$ and $\sigma$ are bounded smooth functions with $||b^{\prime}||_{\infty}<\infty$, $||\sigma^{\prime}||_{\infty} < \infty$. and $||\sigma^{\prime \prime}||_{\infty} < \infty$. Let us  begin first by the simple case when $\sigma$ is constant; that is the additive case.
\subsection{The additive noise case}
We can state the following result
\begin{thm} \label{thm3} Let $\alpha <1$ and $\beta<1$ be  such that $|\rho|<1$ and $\alpha^2 + \beta^2 < \frac{3-2\sqrt{2}}{12}$.  Assume that $b$  is bounded smooth with $||b^{\prime}||_{\infty}< \infty$ and $\sigma(\cdot)\equiv \sigma$ is constant.  Then,  the law of $X_t$ admits a smooth density for all $t\in (0, t_0)$,  where $t_0$ is given by
\begin{equation}\label{eq:t0}
t_0:= \frac{1}{||b^{\prime}||_{\infty}^2} \left(  \frac{(\sqrt{2}-1)^2}{3}- 4(\alpha^2 +\beta^2)\right).
\end{equation}
%\begin{align}\label{smooth:e1}
%C(t_0, \alpha, \beta, b):=3( ||b^{\prime}||^2_{\infty} (t_2-t_1) + 4(\alpha^2 +\beta^2))+2 \sqrt{3( ||b^{\prime}||^2_{\infty}t_0 + 4(\alpha^2 +\beta^2))}<1.
%\end{align}
% Then, for all $t\in (0, t_0]$, the law of $X_t$ admits a smooth density.
\end{thm}
\begin{rk}\label{rk:t0}
Not that thanks to  $\alpha^2 + \beta^2 < \frac{3-2\sqrt{2}}{12}$, we have $t_0>0$. Furthermore, as it is explained in the proof of Theorem \ref{thm3}, the choice of $t_0$ is made possible by ensuring the condition $C(t_0, \alpha, \beta, b)<1$, where the function $C$ is given below by  (\ref{smooth:e1}).  In fact, it is not difficult to see that the  previous condition is equivalent to 
\[\sqrt{3( ||b^{\prime}||^2_{\infty} t_0 + 4(\alpha^2 +\beta^2))}<\sqrt{2}-1,\]
and, clearly, by setting $t_0$ close to zero  the possible values of $\alpha^2 + \beta^2$ are restricted to be also small enough. Namely, for parameters $\alpha$ and $\beta$ such that $\alpha^2 + \beta^2<\frac{(\sqrt{2}-1)^2}{12}$, one can choose such $t_0$. This also explains why we have restricted the range of values of $\alpha$ and $\beta$ in the statement of Theorem \ref{thm3}. 
%
%OR
%
%Not that thanks to  $\alpha^2 + \beta^2 < \frac{3-2\sqrt{2}}{12}$, we have $t_0>0$. Furthermore, as it is explained in the proof of Theorem \ref{thm3}, the choice of $t_0$ is made possible by  considering the function $C$ given by  (\ref{smooth:e1}) and ensuring the condition $C(t_0, \alpha, \beta, b)<1$. In fact, it is not difficult to see that the  previous condition is equivalent to 
%\[\sqrt{3( ||b^{\prime}||^2_{\infty} t_0 + 4(\alpha^2 +\beta^2))}<\sqrt{2}-1,\]
%and, clearly, by setting $t_0$ close to zero  the possible values of $\alpha^2 + \beta^2$ are restricted to be also small enough. Namely, for parameters $\alpha$ and $\beta$ such that $\alpha^2 + \beta^2<\frac{(\sqrt{2}-1)^2}{12}$, one can choose such $t_0$. This also explains why we restrict the range of values of $\alpha$ and $\beta$. 
\end{rk}

For the proof of this result we will check the condition of  Corollary \ref{cor:1}.  We will need the following lemma.
\begin{lem}   \label{lem1: smooth} Let  $b$ be a bounded smooth function such that $|| b^{\prime}||_{\infty}< \infty$ and define
\begin{align}\label{smooth:e1}
C(t, \alpha, \beta, b):=3( ||b^{\prime}||^2_{\infty} t + 4(\alpha^2 +\beta^2))+2 \sqrt{3( ||b^{\prime}||^2_{\infty}t + 4(\alpha^2 +\beta^2))},
\end{align}
for all $t>0$.Then, we have a.s.
\begin{itemize}
\item[(i)] for all $0<t_1<t_2\leq t$, 
\begin{equation} \label{smooth:e2}
| ||DX_{t_2}||_{H}^2-||DX_{t_1}||_{H}^2 | \leq  C(t_2-t_1, \alpha, \beta, b)\times\sup_{0\leq s\leq t}||DX_s ||_{H}^2.
\end{equation}
\item[(ii)] \begin{equation}\label{smooth:e3}
\sup_{0\leq s\leq t}||DX_s||^2_{H} \geq \frac{\sigma^2  t}{2(1+ 3(t^2 ||b^{\prime}||_{\infty}^2 + \alpha^2 +\beta^2))}.
\end{equation} 
\item[(iii)] for all $0 <s  \leq  t$,
 \begin{equation}\label{smooth:e4}
||DX_s ||_{H}^2 \geq \dfrac{1-C(t, \alpha, \beta, b)}{2(1+  3 ||b^{\prime}||^2_{\infty} s^2 +3 (\alpha^2 +\beta^2))}\sigma^2s^2.
\end{equation}
\end{itemize}
\end{lem}
\noindent {\bf Proof of Lemma $\ref{lem1: smooth}$}: 
We have 
\begin{align*}
&| ||DX_{t_2}||_{H}^2-||DX_{t_1}||_{H}^2 |=\left| \int_{0}^{t_2} (D_rX_{t_2})^2dr-\int_{0}^{t_1}(D_rX_{t_1})^2 dr \right|\\
&\leq \int_{t_1}^{t_2} (D_rX_{t_2})^2 dr +  \int_{0}^{t_1} | (D_rX_{t_2})^2-D_rX_{t_1})^2|dr=: J_1+ J_2.
\end{align*}
Let us estimate $J_1$ and $J_2$. For $J_1$, observe that
\begin{align}\label{esmate: J1}
J_1=  \int_{t_1}^{t_2} (D_r X_{t_2}-D_r X_{t_1})^2 dr \leq \int_{0}^{t_2} (D_r X_{t_2}-D_r X_{t_1})^2 dr.
\end{align} 
For $J_2$, using Cauchy-Schwarz inequality and the fact that $t_1 \leq t_2$ , we have 
\begin{align}\label{esmate: J2}
& J_2\leq    \left(\int_{0}^{t_1} (D_rX_{t_2}-D_rX_{t_1})^2dr\right)^{1/2}\left(\int_{0}^{t_1} (D_rX_{t_2}+D_rX_{t_1})^2dr\right)^{1/2}\nonumber\\
&\leq \sqrt{2}\left(\int_{0}^{t_1} (D_rX_{t_2}-D_rX_{t_1})^2dr\right)^{1/2}\left(\int_{0}^{t_1} (D_rX_{t_2})^2+(D_rX_{t_1})^2dr\right)^{1/2}\nonumber\\
&\leq  \sqrt{2}\left(\int_{0}^{t_1} (D_rX_{t_2}-D_rX_{t_1})^2dr\right)^{1/2}\left(\int_{0}^{t_2} (D_rX_{t_2})^2+(D_rX_{t_1})^2dr\right)^{1/2}\nonumber\\
&\leq 2 \sup_{0\leq s\leq t} || DX_s||_H  \left(\int_{0}^{t_1} (D_rX_{t_2}-D_rX_{t_1})^2dr\right)^{1/2}.
\end{align} 
So, by (\ref{esmate: J1}) and (\ref{esmate: J2}), it remains to estimate $\int_{0}^{t_1} (D_rX_{t_2}-D_rX_{t_1})^2dr$. On one hand,  for $r\leq t$, we write 
\begin{align}\label{derivationEquation}
&D_rX_t =\sigma + \int_{0}^{t} D_rb(X_s)ds + \alpha D_r(\max_{0\leq s \leq t}X_s) + \beta D_r(\min_{0\leq s \leq t}X_s).
\end{align}
Hence
\begin{align}\label{estimation: smoothness}
 &(D_rX_{t_2}-D_rX_{t_1})^2\leq  3 |\int_{t_1}^{t_2} D_rb(X_s)ds |^2 + 3 \alpha^2 |D_r(\max_{0\leq s \leq t_1}X_s)-D_r(\max_{0\leq s \leq t_2}X_s)|^2\nonumber\\
&+3 \beta^2 |D_r(\min_{0\leq s \leq t_1}X_s)-D_r(\min_{0\leq s \leq t_2}X_s)|^2\nonumber\\
&\leq 3 |\int_{t_1}^{t_2} D_rb(X_s)ds |^2 + 6 \alpha^2 |D_r(\max_{0\leq s \leq t_1}X_s)|^2+  6 \alpha^2|D_r(\max_{0\leq s \leq t_2}X_s)|^2\nonumber\\
&+6 \beta^2 |D_r(\min_{0\leq s \leq t_1}X_s)|^2+ 6 \beta^2|D_r(\min_{0\leq s \leq t_2}X_s)|^2.
\end{align}
On the other hand, using Cauchy-Schwartz inequality, $(\ref{derivation:b})$ and Fubini theorem, we have
\begin{align}\label{estim:e11}
\int_{0}^{t_2} \left|\int_{t_1}^{t_2} D_rb(X_s)ds\right|^2 dr \leq ||b^{\prime}||^2_{\infty} (t_2-t_1)\sup_{0\leq s\leq t}||DX_s||^2_{H}.
\end{align}
Furthermore, using Proposition \ref{prop1:Nualart} and the fact that $t_1 \leq t_2 \leq t$,  we have 
 \begin{align*}
  \int_{0}^{t_2} \left| D_r(\max_{0\leq s \leq t_1}X_s )\right|^2 dr= || D_r(\max_{0\leq s \leq t_1}X_s)||^2_H \leq \max_{0\leq s \leq t_1}|| DX_s||^2_H \leq  \max_{0\leq s \leq t}|| DX_s||^2_H,
\end{align*}
\begin{align*}
 \int_{0}^{t_2} \left| D_r(\max_{0\leq s \leq t_2}X_s )\right|^2 dr= || D_r(\max_{0\leq s \leq t_2}X_s)||^2_H \leq \max_{0\leq s \leq t_2}|| DX_s||^2_H \leq  \max_{0\leq s \leq t}|| DX_s||^2_H, 
 \end{align*}
and likewise,
\begin{align*}
 \int_{0}^{t_2} \left| D_r(\min_{0\leq s \leq t_2}X_s )\right|^2 dr= || D_r(\min_{0\leq s \leq t_2}X_s)||^2_H \leq \max_{0\leq s \leq t_2}|| DX_s||^2_H \leq  \max_{0\leq s \leq t}|| DX_s||^2_H, 
 \end{align*}
  \begin{align*}
  \int_{0}^{t_2} \left| D_r(\min_{0\leq s \leq t_1}X_s )\right|^2 dr= || D_r(\min_{0\leq s \leq t_1}X_s)||^2_H \leq \max_{0\leq s \leq t_1}|| DX_s||^2_H \leq  \max_{0\leq s \leq t}|| DX_s||^2_H.
\end{align*}
Therefore integrating $(\ref{estimation: smoothness})$ over $[0, t_2]$ and using the above fifth inequalities, we deduce
\begin{align}\label{estimation: smoothness2}
\int_{0}^{t_1} (D_rX_{t_2}-D_rX_{t_1})^2dr \leq  3( ||b^{\prime}||^2_{\infty} (t_2-t_1) + 4(\alpha^2 +\beta^2))\sup_{0\leq s\leq t}||DX_s||^2_{H}.
\end{align}
Consequently  in view of $(\ref{esmate: J1})$ and $(\ref{esmate: J2})$, we obtain
\[J_1 +J_2\leq \left[3( ||b^{\prime}||^2_{\infty} (t_2-t_1) + 4(\alpha^2 +\beta^2))+2 \sqrt{3( ||b^{\prime}||^2_{\infty} (t_2-t_1) + 4(\alpha^2 +\beta^2))} \right]\sup_{0\leq s\leq t}||DX_s||^2_{H}\]
and then  $(\ref{smooth:e2})$ is established.

For $(ii)$, using the inequality $(a+b)^2 \geq \frac{1}{2}a^2 -b^2$, $(\ref{derivationEquation})$, $(\ref{estim:e11})$ and  (\ref{nom_max=e1}) of Proposition \ref{prop1:Nualart}, we have
\begin{align*}
& ||DX_t||^2=\int_{0}^t(D_rX_t)^2  dr \geq \frac{t}{2} \sigma^2 - \dint_{0}^t \left(\int_{0}^{t} D_rb(X_s)ds + \alpha D_r(\max_{0\leq s \leq t}X_s) + \beta D_r(\min_{0\leq s \leq t}X_s)\right)^2 dr\\
&\geq  \frac{\sigma^2  t}{2}   -3\dint_{0}^{t} \left(\int_{0}^{t} D_rb(X_s)ds \right)^2 dr-3 \alpha^2 \dint_{0}^{t} \left(D_r(\max_{0\leq s \leq t}X_s)  \right)^2dr  -3 \beta^2\dint_{0}^{t}  \left(D_r(\min_{0\leq s \leq t}X_s)  \right)^2 dr\\
&\geq \frac{\sigma^2  t}{2}  -3 t^2 ||b^{\prime}||_{\infty}^2 \sup_{0\leq s\leq t}||DX_s||^2_{H} - 3 \alpha^2 \sup_{0\leq s\leq t}||DX_s||^2_{H} -3 \beta^2 \sup_{0\leq s\leq t}||DX_s||^2_{H}\\
&=\frac{\sigma^2  t}{2} -3 (t^2 ||b^{\prime}||_{\infty}^2 + \alpha^2 +\beta^2)\sup_{0\leq s\leq t}||DX_s||^2_{H} .
\end{align*}
Therefore
\begin{align*}
& \sup_{0\leq s\leq t}||DX_s||^2_{H} \geq \frac{\sigma^2  t}{2(1+ 3(t^2 ||b^{\prime}||_{\infty}^2 + \alpha^2 +\beta^2))}.
\end{align*}
For the last claim $(iii)$. Let  $s \in (0, t]$. In view of $(\ref{smooth:e2})$, we have 
\begin{align*}
& ||DX_s||^2_{H} \leq |||DX_s||^2_{H}-||DX_t||^2_{H}| + ||DX_t||^2_{H} \leq C(t-s, \alpha, \beta, b)\sup_{0\leq s\leq t}||DX_s||^2_{H} + ||DX_t||^2_{H}
\end{align*}
Hence\begin{align*}
&  ||DX_t||^2_{H} \geq (1-C(t, \alpha, \beta, b))\sup_{0\leq s\leq t}||DX_s||^2_{H},
\end{align*}
since the function $s\rightarrow C(s, \alpha, \beta, b)$ is increasing. 
Taking into account $(\ref{smooth:e3})$ we deduce $(\ref{smooth:e4})$. This finishes the proof of Lemma \ref{lem1: smooth}.\\

\noindent {\bf Proof of Theorem \ref{thm3}}:  Let $t_0$ be given by (\ref{eq:t0}). Then, according to Remark \ref{rk:t0},  we have  $C(t_0, \alpha, \beta, b)<1$. From Proposition \ref{generalization of thm1} we know that $X_t \in \mathbb{D}^{k, 2}$ for any $t> 0$ and integer $k\geq 1$. Furthermore,  by the point $(iii)$ of Lemma \ref{lem1: smooth} we have a.s., for all $t\in (0, t_0)$,
\begin{equation*}
||DX_t ||_{H}^2 \geq \dfrac{1-C(t_0, \alpha, \beta, b)}{2(1+  3 ||b^{\prime}||^2_{\infty} t^2 +3 (\alpha^2 +\beta^2))}\sigma^2t^2>0.
\end{equation*}
Hence we can apply Corollary \ref{cor:1}; consequently  Theorem \ref{thm3} is proven.
\begin{rk} The statement of Theorem  \ref{thm3} still valid when the function  $b$ is only Lipschitz continuous with Lipschitz constant $K>0$,  by choosing $t_0$ to be
\[t_0:= \frac{1}{K^2} \left( \frac{(\sqrt{2}-1)^2}{3}- 4(\alpha^2 +\beta^2)\right),\]
instead of (\ref{eq:t0}).
\end{rk}
\subsection{The multiplicative noise case}
In this subsection we assume that $\sigma $ is bounded smooth with  $||\sigma^{\prime}||_{\infty}< \infty$, $||\sigma^{\prime \prime}||_{\infty}< \infty$ and  $\inf_{x\in \mathbb{R}}|\sigma(x)|>0$.  Whence either $\inf_{x\in \mathbb{R}}\sigma(x)>0$ or $\sup_{x\in \mathbb{R}}\sigma(x)<0$. Without loss of generality we suppose that
$$\inf_{x\in \mathbb{R}}\sigma(x)>0.$$
In order to state the result of this subsection,  we introduce the following transformation 
\begin{align}\label{transfaormation: e1}
G(y):= \int_{x}^{y} \dfrac{1}{\sigma(u)} du, \,\,\, \mbox{for} \, \, y\in \mathbb{R}.
\end{align}
Clearly $G$ is strictly increasing function with bounded derivative. Hence
\begin{align}\label{max and min}
 \sup_{0\leq s \leq t}G(X_s)= G(\sup_{0\leq s \leq t}X_s)\, \, \mbox{and}\, \,  \inf_{0\leq s \leq t}G(X_s)= G(\inf_{0\leq s \leq t}X_s).
 \end{align}
We also consider the function 
\begin{align}\label{transfaormation: e2}
\tilde{b}(x):= \dfrac{b(G^{-1}(x))}{\sigma(G^{-1}(x))} - \frac{1}{2} \sigma^{\prime}(G^{-1}(x)), \,\,\, \mbox{for} \, \, x\in \mathbb{R}.
\end{align}
In is easily seen that the function $\tilde{b}$ is bounded smooth with $||\tilde{b}^{\prime}||_{\infty}<\infty$.

We can now prove the following theorem
\begin{thm} \label{thm4} Let $\alpha <1$ and $\beta<1$ such that $|\rho|<1$ and $\alpha^2 + \beta^2 < \frac{3-2\sqrt{2}}{12}$.  Let $C$, $G$ and $\tilde{b}$ be the functions introduced in   $(\ref{smooth:e1})$,  $(\ref{transfaormation: e1})$, $(\ref{transfaormation: e2})$ respectively. Then, the law of $X_t$ admits a smooth density for all $t\in (0, t_0)$, where $t_0$ is given by (\ref{eq:t0}).
\end{thm} 
\noindent {\bf Proof of Theorem \ref{thm4}}: By It\^{o}'s formula \cite{RY1999}, we have for all $t\in (0, t_0)$
\begin{align} \label{Ito}
G(X_t)&=  W_t  + \dint_{0}^{t} \left(  \frac{b(X_s)}{\sigma(X_s)}-\frac{1}{2} \sigma^{\prime}(X_s)\right) ds + \alpha \dint_{0}^{t}\frac{1}{\sigma(X_s)}d(\max_{0\leq u\leq s}X_u) + \beta \dint_{0}^{t}\frac{1}{\sigma(X_s)}d(\min_{0\leq u\leq s}X_u)\nonumber\\
&= W_t + \dint_{0}^{t} \tilde{b}(G(X_s))ds + \alpha \dint_{0}^{\max_{0\leq u\leq t}X_u}  \frac{1}{\sigma(u)}du + \beta \dint_{0}^{\min_{0\leq u\leq t}X_u} \frac{1}{\sigma(u)}du\nonumber\\
&= W_t + \dint_{0}^{t} \tilde{b}(G(X_s))ds + \alpha G(\max_{0\leq u\leq t}X_u) +  \beta G(\min_{0\leq u\leq s}X_u)\nonumber\\
&= W_t + \dint_{0}^{t} \tilde{b}(G(X_s))ds + \alpha \max_{0\leq u\leq t}G(X_u) +  \beta\min_{0\leq u\leq t} G(X_u),
\end{align}
where we have used successively the fact that the measures  $d(\max_{0\leq u\leq \cdot}X_u)$ and \\  $d(\min_{0\leq u\leq \cdot}X_u)$  are supported respectively by the sets $\{ s \in (0, t]: X_s =\max_{0\leq u \leq s}X_u\}$  and $\{ s \in (0, t]: X_s =\min_{0\leq u \leq s}X_u\}$, change of variable formula and  $(\ref{max and min})$.
Now applying  Theorem \ref{thm3} of the additive noise case to the doubly perturbed $(\ref{Ito})$, the law of  $G(X_t)$  admits a smooth density and hereby  does the law of $X_t$ since $G$ in invertible. This finishes the proof of Theorem \ref{thm4}.

\end{document}